\documentclass{article}
\usepackage[utf8]{inputenc}
\usepackage{graphicx}
\usepackage{amssymb}
 \usepackage{amsthm}
 \usepackage{amsmath}
 \usepackage{color}

\usepackage{amsfonts}
\usepackage{amssymb}
\usepackage{mathtools}
\usepackage{cite}
\usepackage{hyperref}
\usepackage{enumitem}
\usepackage{blindtext}
\DeclareUnicodeCharacter{00A0}{ }
\DeclareUnicodeCharacter{00A0}{~}
\newtheorem{assumption}{Assumption}

\newcommand{\y}{y}
\newcommand{\s}{s}
\newcommand{\f}{f}

\newtheorem{theorem}{Theorem}

\newtheorem{corollary}[theorem]{Corollary}

\newtheorem{lemma}[theorem]{Lemma}

\newtheorem{proposition}[theorem]{Proposition}

\everymath{\displaystyle}

\providecommand{\U}[1]{\protect\rule{.1in}{.1in}}
\providecommand{\U}[1]{\protect\rule{.1in}{.1in}}

\everymath{\displaystyle}
\usepackage{graphics,color}

\newcommand\gr[1]{{\color{black}{#1}}}

\title{Nonlinear wave interactions \\ in geochemical modeling.}

\author{A. C. Alvarez, J. Bruining, D. Marchesin} 

\date

\begin{document}

\maketitle

\providecommand{\keywords}[1]
{
  \small	
  \textbf{\textit{Keywords---}} #1
}

\begin{abstract}

This paper is concerned with the study of the main wave interactions in a system of conservation laws in geochemical modeling. 
We study the modeling of the chemical complexes on the rock surface. The presence of stable surface complexes affects the relative permeability. We add terms representing surface complexes to the accumulation function in the model presented in \cite{lambert2019nonlinear1}. This addition allows to take into account the interaction of ions with the rock surface in the modeling of the oil recovery by the injection of carbonated water. Compatibility hypotheses with the modeling are made on the coefficients of the system to obtain meaningful solutions.
We developed a Riemann solver taking into account the complexity of the interactions and bifurcations of nonlinear waves. Such bifurcations occur at the inflection and resonance surfaces. We present the solution of a generalized eigenvalue problem in a (n+1)-dimensional case, which allows the construction of rarefaction curves. A method to find the discontinuous solutions is also presented. We find the solution path for some  examples.


 \end{abstract}

\section{Introduction}

Composition modeling  deals with the transport of ions and chemical species in porous media.
Modeling of such processes is based on a system of non-linear partial differential equations composed of Darcy’s law and mass conservation equations. Several processes can be studied with this system aimed at understanding and improving oil recovery (see \cite{lake2002geochemistry}). In this sense the study of the Riemann problem for such flows plays a major role in the field of Mathematics as well as in Engineering.

The class of models described in this paper has first appeared in problems of geochemical flow in porous media, see
 for example \cite{blom2016low,alvarez2016riemann}. In the present work we focus on the mathematical features of the above mentioned physical model. This system represents conservation laws for a set of $n+1$ independent chemical species.
 These chemical species depend on a set of $m$ ion compound concentrations,  which are assumed to have reached thermodynamic equilibrium. This situation leads to a system of $n+1$ nonlinear 
partial equations and $k$ algebraic constraint equations. In this paper, we study the Riemann problem
for this kind of system. 

We focus on mathematical features of systems of conservation laws in the manifold given by the constraints. To simplify the above problem and to obtain an equivalent system with only $n+1$ conservation laws, a novel strategy based on Gibbs rule that reduces  the dimension of the system is developed in \cite{blom2016low,lambert2019nonlinear1}. This procedure was used successfully in studying the injection of low salinity water.

This situation reduces to a system of conservation laws of the type
\begin{equation}
 \displaystyle{\frac{\partial \textbf{G}(U)}{\partial t}+\frac{\partial u\textbf{F}(U)}{\partial x}=0,}\label{eqt}
\end{equation}
with state variables $U=(U_1,\ldots,U_{n})$ belonging to phase space $\Omega$, with $U_i(x,t) : \mathbb{R} \times [0, + \infty) \longrightarrow  \mathbb{R}$, accumulation functions $\textbf{G}(U)=(G_1(U),\cdots,G_{n+1}(U))^T:\Omega \longrightarrow \mathbb{R}^{n+1}$, flux functions
$\textbf{F}(U)=(F_1(U),\cdots,F_{n+1}(U))^T:\Omega \longrightarrow \mathbb{R}^{n+1}$ and  Darcy velocity $u=u(x,t):\mathbb{R}\times [0,+ \infty) \longrightarrow [0,+ \infty)$. The physical and chemical details of this model are described in \cite{alvarez2018analytical1}, while in this work we focus on mathematical and numerical aspects.

System \eqref{eqt} is typical of two phase multicomponent  models for 
geochemical one-dimensional flow, see,  e.g.,  \cite{pope81}. 
In this paper, we focus on the system with accumulation and flux functions given by
\begin{equation}
 G_j=\varphi \rho_{w j}s_w+ \varphi \rho_{o j}s_o+(1-\varphi)\rho_{r j}\quad
 \text{ and }
 \quad 
  F_j=\rho_{w j}f_w+\rho_{o j}f_o, \label{fluxac} \quad j=1,\cdots,n+1
\end{equation}  
where $\varphi$ is a given parameter, $s_w+s_o=1$ and $f_w+f_o=1$. In order to simplify 
the notation we use $s=s_w$, $f=f_w$, 
$s_o=1-\s$ and $f_o=1-f$. 
 
In such a model the state variables are the composition variables $\y_i$, for $i=1,\cdots,n-1$, the saturation
variable  $s_{w}$ and the Darcy velocity $u$. The coefficients $\rho_{w j}$, $\rho_{o j}$ and $\rho_{r j}$ depend on the compositions $\y=(\y_1,\cdots,\y_{n-1})$.
Here $f_\alpha$, which is called
the \textit{fractional flow} of phase $\alpha$, for $\alpha=w,o$, depends on $s$
and $\y=(\y_1,\cdots,\y_{n-1})$ and are written as
$f_\alpha=\lambda_\alpha/(\lambda_w+\lambda_o)$, for $\alpha=w$ and $o$, 
in which {$\lambda_\alpha=k_{r\alpha}/\mu_\alpha$}. The functions 
$\lambda_\alpha$, $k_{r\alpha}$ and $\mu_\alpha$ 
are  called mobility, relative permeability
and viscosity of phase $\alpha$ which depend on the variables $s$ and $y$. The subscript $w$ in 
$\rho_{wj}$ indicates that these functions refer to the water phase while the subscript $o$ for $\rho_{oj}$ indicates the oil phase. Finally, $\rho_{rj}$ represent the ion interactions with the porous medium.

 Here, we incorporate to the model presented in \cite{lambert2019nonlinear1}, the coefficients $\rho_{rj}$ in the accumulation functions $G$, which arise from the interaction of ions with the rock surfaces. Such chemical interaction enables incorporating the formation of surface complexes, which enhances the oil recovery by injection of low salinity carbonated water. Therefore, when these coefficients are introduced new features appear in the Riemann  solutions associated to the system \eqref{eqt}.

In our study we consider the physical state variables $\Omega=\{(\s,\y)= [0,1]\times\mathcal{K}\subset \mathbb{R}\times \mathbb{R}^{n-1}\}$, where $\mathcal{K}$ denotes the $n-1$-dimensional hypercube. The physical model consists of $n-1$ species concentrations and two additional variables $(s,u)$, leading to $n+1$ unknowns $U=(s,y,u)$.

A Riemann problem consists of a Cauchy (initial value) problem governed by equations of type \eqref{eqt} with initial data 
\begin{equation}
U(x,t=0) = \left \{ \begin{matrix} U_L & \mbox{if } x < 0,
\\ U_R & \mbox{if } x > 0.\end{matrix}\right.
\label{eqrt2}
\end{equation}

 Riemann solutions are obtained using rarefaction, shock and composite wave curves going from the state $U_L$ to $U_R$ 
 passing through a set of intermediate constant states. Rarefactions are continuous self-similar solutions of \eqref{eqt}, which are represented by
\begin{equation}
U=\widehat{U}(\xi),\quad \text{with} \quad \xi=x/t.
\label{eqt2}
\end{equation}
Substituting \eqref{eqt2} into system \eqref{eqt} we obtain the generalized eigenvalue problem
\begin{equation}
Ar=\lambda Br, \quad \text{where} \quad A={\partial F}/{\partial U}, \quad B={\partial G}/{\partial U}.
\label{ge1}
\end{equation}
The eigenvector $r$ is parallel to $d\widehat{U}/d\xi$, so the rarefaction curves are tangent to the characteristic field given by the normalized eigenvector $r$.  
 
A shock wave is a traveling discontinuity in a (weak) solution of system \eqref{eqt} given by
\begin{equation}
U(x,t) = \left \{ \begin{matrix} U^- & \mbox{if } x < \sigma t,
\\ U^+ & \mbox{if } x > \sigma t,\end{matrix}\right.
\label{eqrt2a}
\end{equation}
where $\sigma$ is a real constant called the shock speed. Solution \eqref{eqrt2a} is a piecewise constant weak solution for the Riemann problem  defined by \eqref{eqt} and \eqref{eqrt2} provided $U^+$ and $U^-$ satisfy the Rankine-Hugoniot
condition:
\begin{equation}
F (U^-) - F (U^+) = \sigma(G(U^-) - G(U^+)).
\end{equation}
In this work we use the wave curve method (see \cite{liu1974riemann,liu1975riemann,azevedo1995multiple}), 
a valuable tool to obtain semianalytical solutions.
More details about the wave curve method can be found in \cite{lax1957hyperbolic,azevedo1995multiple,lambert2010riemann,dahmen2005riemann}, with theoretical justification provided by  \cite{liu1974riemann,liu1975riemann,issacson1992global,wendroff1972riemann,alvarez2020resonance} and references therein. Here, we provide a method for calculating and concatenating waves in a set of relevant examples. To do this, we provide the waves and the bifurcation surfaces necessary to solve the Riemann problem for the system of conservation laws \eqref{eqt}. To study this system in the (n+1)-dimensional case is a hard task in general, however due to the structure of the flux and accumulation functions we can obtain some useful information about rarefaction and shock waves. 

In general the system \eqref{eqt} is not strictly hyperbolic, i.e., the generalized eigenvalues in \eqref{ge1} must coincide on some hypersurface. Moreover, genuine nonlinearity fails, i.e., there is a locus where the characteristic speed has zero directional derivative along the vector field $r=r_i(U)$, i.e.  $\lambda^{\prime}_i(r)= \nabla \lambda_i \cdot r_i$ = 0. The rarefaction curve generically stops at this locus. Thus the determination of the inflection locus is relevant for the wave curve method.

Following Liu (see \cite{liu1974riemann,liu1975riemann}) to cross the inflection locus, the rarefactions wave need to be attached to a shock. This corresponds to the construction of a composite curve, which arises in state space during the construction of wave curves for the solution of a Riemann problem for  non-strictly hyperbolic systems of conservation laws. 

There is one wave in which the saturation changes (called the Buckley-Leverett wave (B-L)). There are also $n-1$ chemical waves.  Having identified these waves we provide the analysis of the Riemann solutions in certain relevant situations.

An in-house Riemann solver was developed to represent wave curves that satisfy compatibility and admissibility criteria (see \cite{alvarez2020resonance} for technical details of implementation). In this solver, we automate the construction of the solution paths taking into account the bifurcation surfaces.

The paper is organized as follows. Section \ref{fraction} presents the fractional flow function used here. In Section $\ref{secteog}$, we obtain the characteristic speeds and characteristic vectors of the system. In Section \ref{rp} we obtain the bifurcation surfaces crucial in determining the wave interactions. In Section $\ref{rsk}$, we obtain the Hugoniot locus. This locus is the main ingredient to obtain 
the Riemann solutions. In Section \ref{RP1}, we solve the  Riemann  problem for some examples. 
In two numerical examples, we study the case when the coefficients $\rho_{jr}$ are constant. In  Section \ref{con}, we draw our conclusions.

\section{Fractional flow functions}
\label{fraction}
 
  The fractional flows $f_w$ and $f_o$ for water and oil in \eqref{fluxac} are frequently taken as saturation-dependent functions defined as follows. We denote
 $s_e=(s-s_{wc})/(1-s_{wc}),~\text{for}~
 s \geq s_{wc}~\text{and}~s_e=0~ \text{for}$ $s < s_{wc}; k_{rw}(s)=s_e^{(2/\lambda+3)}~~ \text{ and }~~ k_{ro}(s_o)=(1-s_e)^2(1-s_e^{(2/\lambda+1)}),$ ($s_o=1-s$) (\cite{brooks1963hydraulic}). Here, the water viscosity is taken as $\mu_w=0.001$ and the oil viscosity
 as $\mu_o=0.002$ when they are constant. In other examples $\mu_w$ and $\mu_o$ depend on $y$ also. The fractional flow functions for the water and oil phases are given by
 \begin{equation}
 f_w(s)=\frac{k_{rw}(s)/\mu_w}{(k_{rw}(s)/\mu_w+k_{ro}(s_o)/\mu_o)},\quad \text{and} \quad f_o(s)=1-f_w(s).
 \label{frac1ab}
 \end{equation}
The fractional flow $f_w$ depends on the chemical variables via the oil viscosity. In order to take into account this dependence in the analysis of the wave interaction, we assume that the fractional flow
$f=f_w$ is $S$-shaped (see \cite{buckley1942mechanism}):
\begin{assumption}
 \label{assumption1}\textbf{S-shape}.
 We assume that  the fractional
flows  $\f(\s,\y)$ has the following suitable properties:
$\f(0,y)=0$ and $\f(1,y)=1$, $\frac{\partial \f(0,y)}{\partial \s}=\frac{\partial \f(1,y)}{\partial \s}=0$.
For each $y$, $f$ has one inflection point $s^*$ (that depends on $\y$ and possibly on viscosities)
\begin{equation}
s^*=s^*(\y)\quad \text{such that} \quad
\frac{\partial^2 \f(s^*,\y)}{\partial \s^2}=0.\label{star}
\end{equation}
For each $y$, there are two states 
\begin{equation}
s^{\dag,1}(\y)<s^*(\y)<s^{\dag,2}(\y)\quad  \text{such that}
\quad \frac{\partial \f(s^{\dag,1},\cdot)}{\partial \s}=
\frac{\partial \f(s^{\dag,2},\cdot)}{\partial \s}=1.\label{s1s2}
\end{equation}
For  $0<s<1$, $\displaystyle{\frac{\partial \f(s,\cdot)}{\partial \s}> 0}$. If $0\leq s <s^*(\y)$, then 
 $\displaystyle{\frac{\partial^2 \f(s,\y)}{\partial \s^2}>0}$, and if  $s^*(\y)< s \leq 1$, then    $\displaystyle{\frac{\partial^2 \f(s,\y)}{\partial \s^2}<0}$.
 \label{assum2}
\end{assumption}

\section{Eigenvalues analysis and elementary waves}
\label{secteog}

The system has $n+1$ equations with the unknown variables $ (s, y, u) $
with $ y = (y_1,y_{2},\cdots,y_{n-1})$ given by
\begin{align}
&\frac{\partial }{\partial t}( \varphi \rho_{wj}(y)s_w+\varphi \rho_{oj}(y)s_o+(1-\varphi)\rho_{rj}(y))+
\frac{\partial }{\partial x}(u(\rho_{wj}(y)f_w+\rho_{oj}(y)f_o))=0.\label{eq1}
\end{align}
We take $s=s_w$, $f=f_w$, $s_o=1-s$ and $f_o=1-f_{w}$.
The accumulation and flux functions $G$ and $uF$ are given in terms of
\begin{align}
&G_j=\varphi \rho_{wj}(y)s_w+\varphi \rho_{oj}(y)s_o+(1-\varphi)\rho_{rj}(y)\\
&F_j=\rho_{wj}(y)f_w+\rho_{oj}(y)f_o
\label{acful}
\end{align}
The index $w$ (water) is often replaced by the index $a$ (aqueous phase), and the index $j=1,\cdots,n+1$ is used to denote chemical species.

The main object to obtain the Riemann 
solutions of hyperbolic systems are 
the rarefaction and shock curves.
The rarefactions are obtained from integral
curves of the line fields given by the eigenvectors
of system $A\vec{r}=B\lambda\vec{r}$, with the Jacobian matrices $B=DG$ and $A=D(uF)$ of the accumulation and flux functions $uF$ and $G$. For system $(\ref{eq1})$, the matrices $B$ and $A$ are given by
\begin{align}
&B_{i,1}=\varphi[\rho_i], \quad B_{i,k+1}=\varphi\frac{\partial \rho_{wi}}{\partial y_k}s+\varphi\frac{\partial \rho_{oi}}{\partial y_k}s_o+(1-\varphi)\frac{\partial \rho_{ri}}{\partial y_k},  \text{ and } B_{i,n+1}=0.\label{matrixB} \\
&A_{i,1}= u[\rho_i]\frac{\partial f}{\partial s}, \quad 
A_{i,k+1}= u\left(\frac{\partial \rho_{wi}}{\partial y_k}f
+\frac{\partial \rho_{oi}}{\partial y_k}f_o \right)\text{ and }
A_{i,n+1}= F_i, \label{matrixA}
\end{align}
in which $i=1,\cdots,n+1$, $k=1,n-1$ and
\begin{equation}
[\rho_i]=\rho_{wi}-\rho_{oi}. \label{diferenca}
\end{equation}
From the Jacobian matrices, we obtain the eigenpairs, which we summarize in 
\begin{proposition}
	\label{lemakt}
	The eigenpairs of the eigenvalue problem $A\vec{r}=\lambda B\vec{r}$, where the matrices $B$ and $A$ represent the Jacobian
	of the accumulation and flux terms of system $(\ref{eq1})$  are the Buckley-Leverett  eigenpair ($\lambda_s,\vec{r}_s$)  given by
	\begin{equation}
	\lambda_s=\frac{u}{\varphi}\frac{\partial f}{\partial s}
	\quad \text{ and } \quad \vec{r}_s=(1,\cdots,0)^T,\label{lambdasa}
	\end{equation}
 together with the $n-1$ chemical eigenpairs $(\lambda_{\Lambda_i},\vec{r}_{\Lambda_i})$ given by 
	\begin{equation}
	\lambda_{\Lambda_i}=\frac{u}{\varphi}\frac{f-\Lambda_i}{s-\Lambda_i},\label{eig2a}
	\end{equation}
	with $i=1,\cdots,n-1$.
	We obtain $\Lambda_i$ and $\vec{v}_i$ as the solutions of the  generalized eigenvalue problem 
	\begin{equation}
	((\mathcal{A}+\mathcal{A}_1)-\Lambda_i (\mathcal{B}+\mathcal{B}_1))\vec{v}_i=0,
	\label{eigeqs2}
	\end{equation}
	where the matrices $\mathcal{A}, \mathcal{A}_1, \mathcal{B}$ and  $\mathcal{B}_1$  are given in \eqref{sd1}-\eqref{sd4}. The matrices $\mathcal{A}$ and  $\mathcal{B}$ depend on the variable $y$
	while $\mathcal{A}_1$ and  $\mathcal{B}_1$ depend on the
	variables $s,y$. Moreover, we obtain $\vec{r}_{\Lambda_i}= (r_{\Lambda_i}^1,\cdots,r_{\Lambda_i}^{n+1})^T$ as 
      \begin{equation}
     \vec{r}_{\Lambda_i}=P\vec{v}_i,\\
    \label{eiget1}
    \end{equation}
   where the matrix $P$ is defined in \eqref{eiget1a2}.
    Furthermore, the left eigenvectors $l_{\Lambda_i}$ are given by
    $(0,l_{\vec{v_i}},0)$ where $l_{\vec{v_i}}$ are the left eigenvectors associated to \eqref{eigeqs2}, i.e.,
    	\begin{equation}
	l_{\vec{v_i}}((\mathcal{A}+\mathcal{A}_1)-\Lambda_i (\mathcal{B}+\mathcal{B}_1))=0.
	\label{eigeqs2a}
	\end{equation}

\end{proposition}

\textbf{Proof of Proposition}

The idea of the proof consists of reducing the matrix $\mathcal{G}$ by straightforward application
of the Gauss elimination procedure, i.e., transforming the matrix $\mathcal{G}$ into an equivalent matrix and obtaining the eigenvalues and eigenvectors easily.

To determine the eigenvalues we solve $det(A-\lambda B)=0$,
and to calculate eigenvectors we solve $(A-\lambda B)\vec{r}=0$, where
the matrices $A$ and $B$ are given in \eqref{matrixA}. Let us denote by $\mathcal{G}=(A-\lambda B)$, and $\mathcal{G}=(\mathcal{G}_{i,j})$ for $i,j=1,\cdots,n+1$.
Here we use an auxiliary index $k$ that ranges from $1$ to $n-1$ and we write $\mathcal{G}_{i,j}$ as
\begin{equation}
\mathcal{G}_{i,1}(\lambda)=[\rho_i]\xi_1(\lambda), 
\mathcal{G}_{i,k+1}(\lambda)=\frac{\partial \rho_{wi}}{\partial y_k}\xi_2(\lambda)+
\frac{\partial \rho_{oi}}{\partial y_k}\xi_3(\lambda)- (1-\varphi)\frac{\partial  \rho_{ri}}{\partial y_k}\lambda
\text{ , } 
\mathcal{G}_{i,n+1}= F_i. \label{matrixnew}
\end{equation}
where we define the auxiliary functions
\begin{equation}
\xi_1(\lambda)=\left(u\frac{\partial f}{\partial s}-\varphi \lambda \right), \quad
\xi_2(\lambda)=(u f-\varphi \lambda s)\quad \text{ and }\quad
\xi_3(\lambda)=(uf_o-\varphi \lambda s_o).\label{xis}
\end{equation}
To express calculations done in the  next steps of the proof, it is useful to define the following coefficients $\gamma_{ij}$, 
$\varrho_{ij}$, $\nu_i$, $\vartheta_{ij}$ and $\varsigma_{ij}$  
for $i=1,\cdots,n+1$, $j=1,n-1$ as
\begin{align}
&\gamma_{ij}=\frac{\partial \rho_{wi}}{\partial y_j}[\rho_1]-
\frac{\partial \rho_{w1}}{\partial y_j}[\rho_i],
\quad
\varrho_{ij}=\frac{\partial \rho_{oi}}{\partial y_j}[\rho_1]-
\frac{\partial \rho_{o1}}{\partial y_j}[\rho_i],
\label{coef1}
\end{align}
\begin{align}
&\pi_{ij}=(1-\varphi)(\frac{\partial \rho_{ri}}{\partial y_j}[\rho_1]-
\frac{\partial \rho_{r1}}{\partial y_j}[\rho_i]),
\quad 
\nu_i=[\rho_1]F_{i}-[\rho_i]F_{1},
\label{coef2}
\end{align}
\begin{align}
\label{gamma1}
 &\vartheta_{ij}=\gamma_{ij}\nu_{4}-\gamma_{4,j}\nu_i \quad
 \varsigma_{ij}=\varrho_{ij}\nu_{4}-\varrho_{4,j}\nu_i,\\
&  \tau_{ij}=\pi_{ij}\nu_{4}-\pi_{4,j}\nu_i.
 \label{gamma}
 \end{align}
All these coefficients (\eqref{coef1}- \eqref{gamma}) depend only of the variables $y$. Now, we are able to perform the Gaussian elimination 

\begin{itemize}
	
\item[1)] Substituting  the $i$-th row of matrix $(\ref{matrixnew})$, for
$i=2,\cdots,n+1$,  by the sum of the first row of $(A-\lambda B)$ 
multiplied by $-[\rho_i]$ with its $i$-th row   multiplied by 
$[\rho_1]$, we obtain, for $j=1,n-1$ and $i=2,\cdots,n+1$

\begin{equation}
\begin{pmatrix}
\mathcal{G}_{1,1}(\lambda) & \mathcal{G}_{1,j+1}(\lambda) &\mathcal{G}_{1,n+1}\\
\mathbb{O} & \mathcal{G}^{(1)}_{i,j+1}(\lambda) & \mathcal{G}^{(1)}_{i,n+1}\\
\end{pmatrix}\vec{r}=0.\label{matrix3}
\end{equation}

Here $\mathbb{O}$ is the column vector of three zeros and block matrices $(\mathcal{G}^{(1)}_{i,j+1})$ and $ (\mathcal{G}^{(1)}_{i,n+1})$ for
$i=2,\cdots,n+1$ and $j=1,n-1$ are the block matrices with elements 
\begin{equation}
\mathcal{G}^{(1)}_{i,j+1}(\lambda)=\left(\gamma_{ij}\xi_2(\lambda)+\varrho_{ij}\xi_3(\lambda)-\pi_{ij}\lambda\right)_{1\leq j\leq n-1}\quad 
\text{ and }
\quad 
\mathcal{G}^{(1)}_{i,n+1}=\nu_{i},
\label{eq69}
\end{equation}
where $\gamma_{ij}$, $\varrho_{ij}$ $\pi_{ij}$ and  $\nu_{i}$ are given by Eq. \eqref{coef1} and \eqref{coef2}.

Notice here that if $[\rho_i]$ given in $(\ref{diferenca})$ is zero for some index $i$ at some $y$, the corresponding entry in the first column is zero and we do not need to perform Gaussian elimination. Moreover, if $[\rho_1]$ is zero for some $y$, we interchange the first row  with another row to obtain a non-zero pivot.

\item[2)] Assuming that $\nu_{n+1}$ is non-zero, we substitute the $i$-th row of $(\ref{matrix3})$,i .e., ($\mathcal{G}^{(1)}_{i,j+1}(\lambda)$,
$j=1,\cdots,n$), for $i=2,n$,  by the sum of the $(n+1)$-th row of $(\ref{matrix3})$, ($\mathcal{G}^{(1)}_{n+1,j+1}(\lambda)$,
$j=1,\cdots,n$)
multiplied by $-\nu_i$ with the $i$-th row of $(\ref{matrix3})$ multiplied by $\nu_{n+1}$, 
Notice that the last column from $i=2,\cdots,n$
becomes $\nu_{n+1} (-\nu_i)-\nu_{n+1} \nu_i=0$,
and the first column for $i=2,\cdots,n$ remains as zero. Thus, we obtain for $j=1,\cdots,n-1$ and
$i=2,\cdots,n-1$
\begin{equation}
\begin{pmatrix}
\mathcal{G}_{1,1}(\lambda) & \mathcal{G}_{1,j+1}(\lambda) &\mathcal{G}_{1,n+1}\\
\mathbb{O} & \mathcal{G}^{(2)}_{i,j+1}(\lambda) &\mathbb{O}\\
0 & \mathcal{G}^{(1)}_{n+1,j+1}(\lambda)&\nu_{n+1}\\
\end{pmatrix}\vec{r}=0.\label{pms}
\end{equation}
Here $\mathbb{O}$ is the column vector of $n-1$ zeros.
Also, $\mathcal{G}^{(2)}_{i,j+1}(\lambda)$ and $ \mathcal{G}^{(1)}_{n+1,j+1}(\lambda)$ for
$i=2,\cdots,n$ and $j=1,\cdots,n-1$ are the block matrices given by:
\begin{equation}
\mathcal{G}^{(2)}_{i,j+1}(\lambda)=\vartheta_{ij}\xi_2(\lambda)+\varsigma_{ij}\xi_3(\lambda)-\tau_{ij}\lambda
,~\\
\mathcal{G}^{(1)}_{n+1,j+1}(\lambda)=\gamma_{n+1,j}\xi_2(\lambda)+\varrho_{n+1,j}\xi_3(\lambda)-\pi_{n+1,j}\lambda, \label{g2k}
\end{equation}
where $\vartheta_{ij}$, $\varsigma_{ij}$ and $\tau_{ij}$ are given by \eqref{gamma1} and $(\ref{gamma})$, for $i=1,\cdots,n+1$, $j=1,n-1$. 


\item[3)]
Since $\mathcal{G}$ is a block matrix, it is useful to define the
matrix $\left(G^{(2)}_{l,r}\right)$ for $l,r=1,\cdots,n-1$  from matrix $\left(\mathcal{G}^{(2)}_{i,j+1}\right)$, with rows from $2$ to $n$ and columns from $2$ to $n$.
\begin{equation}
G^{(2)}_{l,r}=\mathcal{G}^{(2)}_{l+1,r+1}(\lambda) \quad \text{ for } l,r=1,\cdots,n-1.\label{gge}
\end{equation}

From $(\ref{pms})$ and $\eqref{gge}$, we have $\det(A-\lambda B)=\mathcal{G}_{1,1}(\lambda)\det({G}^{(2)}_{l,r}) \nu_{n+1}$ 

From $\mathcal{G}_{1,1}(\lambda)=0$ we get $\xi_1(\lambda)=0$, therefore  $(\ref{xis}.a)$ holds.
To determine the corresponding eigenvalue, we substitute $\lambda$ by $\lambda_s$ into $(\ref{pms})$ and  the statement $(\ref{xis}.b)$ holds. Finally, the eigenpair $(\lambda_s,\vec{r}_s)$  given by $(\ref{lambdasa})$ is obtained. 
In this eigenpair, only saturation varies and we call this wave a \textit{saturation wave}, or
Buckley-Leverett type wave. 

The other eigenvalues are obtained by solving
\begin{equation}
\det({G}^{(2)}_{l,r})=0,
\end{equation}
where the matrix ${G}^{(2)}_{l,r}$ is given by \eqref{gge}. 

To determine the other eigenvalues, we substitute $\lambda$ by $\lambda_{\Lambda_{i}}$  
given by $(\ref{eig2a})$ into $(\ref{gge})$. Notice that 
\begin{equation}
\xi_2(\lambda_{\Lambda_i})=u f-\varphi \lambda_{\Lambda_i} s= u\frac{f(s-\Lambda_i)-(f-\Lambda_i)s}{s-\Lambda_i}=u\Lambda_i\frac{s-f}{s-\Lambda_i},\label{ueq1}
\end{equation}
\begin{equation}
\xi_3(\lambda_{\Lambda_i})=uf_o-\varphi \lambda_{\Lambda_i} s_o=u(1-f)-\lambda_{\Lambda_i}(1-s)=u(1-\Lambda_i)\frac{s-f}{s-\Lambda_i}.\label{ueq1n}
\end{equation}

Substituting $\xi_2(\lambda_{\Lambda_i})$ and $\xi_3(\lambda_{\Lambda_i})$ given by \eqref{ueq1}-\eqref{ueq1n}
into \eqref{gge}, the matrix ${G}^{(2)}_{l,r}$ can be rewritten as

\begin{eqnarray*}
{G}^{(2)}_{l,r} = \frac{\lambda_{\Lambda_i}(s-f)}{f-\Lambda_i}\begin{bmatrix} 
\vartheta_{2,1}+\Lambda_i(\vartheta_{2,1}-\varsigma_{2,1})          & \dots         & \vartheta_{2,n-1}+\Lambda_i(\vartheta_{2,n-1}-\,\varsigma_{2,n-1})         \\ 
\vdots   &  \ddots    & \vdots \\
\vartheta_{n,1}+\Lambda_i(\vartheta_{n,1}-\,\varsigma_{n,1})          & \dots         & \vartheta_{n,n-1}+\Lambda_i(\vartheta_{n,n-1}-\,\varsigma_{n,n-1})        
\end{bmatrix}\\
-\frac{\lambda_{\Lambda_i}(s-f)}{f-\Lambda_i} \begin{bmatrix} 
\frac{f-\Lambda_i}{s-f}\tau_{2,1}           & \cdots         & \frac{f-\Lambda_i}{s-f}\tau_{2,n-1}         \\ 
\vdots   &  \ddots    & \vdots \\
 \frac{f-\Lambda_i}{s-f}\tau_{n,1}            & \cdots         & \frac{f-\Lambda_i}{s-f}\tau_{n,n-1}        \\
\end{bmatrix}.\\
\end{eqnarray*}

Let us denote
 \begin{equation}
\mathcal{A}=\begin{bmatrix} 
\vartheta_{2,1}          & \dots         & \vartheta_{2,n-1}        \\ 
\vdots   &  \ddots    & \vdots \\
\vartheta_{n,1}         & \dots         & \vartheta_{n,n-1}      
\end{bmatrix},
\label{sd1}
\end{equation}
 \begin{equation}\mathcal{B}=-\begin{bmatrix} 
(\vartheta_{2,1}-\varsigma_{2,1})          & \dots         & (\vartheta_{2,n-1}-\,\varsigma_{2,n-1})         \\ 
\vdots   &  \ddots    & \vdots \\
(\vartheta_{n,1}-\,\varsigma_{n,1})          & \dots         & \vartheta_{n,n-1}-\,\varsigma_{n,n-1})   \end{bmatrix},
 \label{sd2}
 \end{equation}
 
 \begin{equation}\mathcal{A}_1=-\begin{bmatrix} 
\frac{f}{s-f}\tau_{2,1}           & \cdots         & \frac{f}{s-f}\tau_{2,n-1}         \\ 
\vdots   &  \ddots    & \dots \\
 \frac{f}{s-f}\tau_{n,1}            & \dots         & \frac{f}{s-f}\tau_{n,n-1}        \\
\end{bmatrix}, \label{sd3} \end{equation} 
and
 \begin{equation}\mathcal{B}_1=-\begin{bmatrix} 
\frac{-1}{s-f}\tau_{2,1}           & \cdots         & \frac{-1}{s-f}\tau_{2,n-1}         \\ 
\vdots   &  \ddots    & \dots \\
 \frac{-1}{s-f}\tau_{n,1}            & \dots         & \frac{-1}{s-f}\tau_{n,n-1}        \\
\end{bmatrix}. \label{sd4} 
 \end{equation}
Let us denote  $\varUpsilon=(\lambda_{\Lambda_i}(s-f)/(f-\Lambda_i))$, we have
\begin{equation}
{G}^{(2)}_{l,r}=\varUpsilon((\mathcal{A}-\Lambda_i\mathcal{B})+(\mathcal{A}_1-\Lambda_i\mathcal{B}_1)),
\label{e40}
\end{equation}
where the matrices $\mathcal{A}$ and $\mathcal{B}$ depend only on the variables $y$, while  the matrices $\mathcal{A}_1$ and $\mathcal{B}_1$ depend on $y$ and $s$.

Thus to calculate the eigenvalues $\lambda_{\Lambda_i}$ it is sufficient to determine $\Lambda_i$ using the generalized eigenvalue problem given in \eqref{eigeqs2}.

\item[4)] The right eigenvectors $\vec{r}_{\Lambda_i}$ related to $\lambda_{\Lambda_i}$ are obtained by substituting 
$\lambda$ by $\lambda_{\Lambda_i}$ in \eqref{pms}.
Because of the structure of the matrix in \eqref{pms}, we can split the calculation of $\vec{r}_{\Lambda_i}$. 
First, we obtain 
the coordinates $(r_{\Lambda_i}^2,\ldots,r_{\Lambda_i}^{n})$ of eigenvector 
$\vec{r}_{\Lambda_i}$. 
Let us define the auxiliary vector $\vec{v}_i$
of size $n-1$, which solves
\begin{equation}
\left({G}^{(2)}_{l,r}\right)\vec{v}_i=0.
\label{pmsnok}
\end{equation}
Using  $(\ref{ueq1})$-$(\ref{ueq1n})$, we can see that
$(\ref{pmsnok})$, after simplifications is written in the form $(\ref{eigeqs2})$. Second, the coordinates $r_{\Lambda_i}^{n+1}$ and $r^1_{\Lambda_i}$ are  calculated by
solving the first and the last equations of $(\ref{pms})$ as
\begin{align}
     \label{eiget1a}
    &r_{\Lambda_i}^{n+1}=-\frac{\displaystyle \sum_{j=2}^{n}\mathcal{G}^{(1)}_{n+1,j}(\lambda_{\Lambda_i})v_i^{j-1}}{\mathcal{G}_{n+1,n+1}},\\
    &r^1_{\Lambda_i}=-\frac{\mathcal{G}_{1,n+1}r_{n+1}+\displaystyle \sum_{j=2}^{n}\mathcal{G}_{1,j}(\lambda_{\Lambda_i})v_i^{j-1}}{\mathcal{G}_{1,1}(\lambda_{\Lambda_i})}.	
    \label{eiget1a1}
    \end{align}
    and the coordinates $r_{\Lambda_i}^j$ for $j=2,\cdots,n$ are given by
    \begin{equation}
         r^j_{\Lambda_i}=v^{j-1}_i,
    \end{equation}
    where $v_i$ is given in \eqref{pmsnok}.
    
    Here the matrices $\mathcal{G}$ and $\mathcal{G}^{(1)}_{n+1,j}$ evaluated at $\lambda=\lambda_{\Lambda_i}$ are given in \eqref{matrixnew} and  \eqref{eq69}.

Let us denote $a=(a_1,\cdots,a_{n-1})$ and  $b=(b_1,\cdots,b_{n-1})$ with
\begin{equation}
    a_j=-\mathcal{G}^{(1)}_{n+1,j}(\lambda_{\Lambda_i})/\mathcal{G}_{n+1,n+1},
    \label{eA}
\end{equation}
 and
 \begin{equation}
    b_j=\mathcal{G}_{1,n+1}\mathcal{G}^{(1)}_{n+1,j}(\lambda_{\Lambda_i})/(\mathcal{G}_{n+1,n+1}\mathcal{G}_{1,1})-\mathcal{G}_{1,j}(\lambda_{\Lambda_i})/\mathcal{G}_{1,1}. 
    \label{eB}
 \end{equation}

Let us denote the matrix $P=(P_{i,j})$ with $n+1$ rows and $n-1$ columns as
\begin{equation}
P_{i,j}=\delta_{i,1}a_j+\delta_{i,n+1}b_j+\delta_{i,j+1}, \label{eiget1a2}
\end{equation}
with  $i=1,\ldots,n+1$, $j=1\ldots,n-1$, $\delta_{i,j}=1$ if $i=j$,
and $\delta_{i,j}=0$ otherwise.

Using \eqref{eiget1a2} and \eqref{eiget1a}- \eqref{eiget1a1}
we obtain \eqref{eiget1}.$\square$
\end{itemize}
The integral curves $\mathcal{W}^s$  and $\mathcal{W}^{\Lambda_i}$ associated to $\vec{r}_s$ and $\vec{r}_{\Lambda_i}$ are obtained by integrating the ODEs 
  \begin{equation}
     \frac{d(s,y,u)}{d \xi}=\vec{r}_s,
     \label{ode1}
 \end{equation}
  and 
\begin{equation}
\frac{d(s,y,u)}{d \xi}=\vec{r}_{\Lambda_i}.
\label{integral}
\end{equation}
Finally, we conclude that in the Buckley-Leverett wave only the saturation changes. Moreover, it is possible to verify the following
 \begin{corollary}
 \label{corolinfle}
  If $\lambda_\Lambda=u{f}/{\varphi s}$ ($\Lambda_i=0$) or $\lambda_\Lambda=u{(1-f)}/\varphi (1-s)$ ($\Lambda_i=1$), then
   the field associated to $\lambda_\Lambda$ is linearly degenerate.
 \end{corollary}
 
\section{Bifurcation surfaces}
\label{rp}

\subsection{Resonance}

\textit{Bifurcation loci} are used to divide the phase space in subregions in which the sequences of waves for the Riemann solutions are the same. To study these loci we take the flux functions $f$ given in \ref{fraction}.

For each $\Lambda_i$ there exist a \textit{coincidence locus}, which is denoted by $\Gamma_{\Lambda_i}$; it occurs where the eigenvalues $\lambda_s$ and 
$\lambda_{\Lambda_i}$ coincide, i.e., at the zero of
\begin{equation}
 \mathcal{G}_{\Lambda_i}(s,y)=\frac{\partial  \f(s,y)}{\partial \s} - \left(\frac{\f(s,y)-\Lambda_i(s,\y)}{s-\Lambda_i(s,\y)}\right).\label{gfun}
\end{equation}
A relevant situation occurs when the system of conservation laws in \eqref{eq1} has constant coefficients $\rho_{rj}$ with $j=1,\cdots,n+1$, which we summarize in the corollary
\begin{corollary}
In the system of conservation laws \eqref{eq1},
        if the coefficients $\rho_{rj}$ with $j=1,\cdots,n+1$ are constant then the eigenvalues $\Lambda_i$
        in \eqref{eigeqs2} depend only on the variables $y$.
        \label{yes}
\end{corollary}
 If any coefficient $\rho_{rj}$ are not constant
then there occur significant changes in the behavior of the bifurcation loci. In particular notice that  we cannot provide a characterization of the inflection locus  $\mathcal{G}_{\Lambda_i}(s,y)=0$ in the nonconstant case. However making some assumptions on the coefficients of the system of conservation laws  \eqref{eq1}, we have that each $\Gamma_{\Lambda_i}$ consists of the union of two disconnected surfaces.
\begin{lemma}
  \label{lema1} 
  Assume that the coefficients $\rho_{rj}$ with $j=1,\cdots,n+1$ in \eqref{eq1} are constant. Assume that either $\Lambda_i<0$ or $\Lambda_i>1$. We conclude $\Gamma_{\Lambda_i}=\Gamma_{1,\Lambda_i}\cup\Gamma_{2,\Lambda_i}$, with
  \begin{align}
  \Gamma_{1,\Lambda_i}=\{(s_{\Lambda_i}^1,y) : \mathcal{G}_{\Lambda_i}(s_{\Lambda_i}^1(\y),y)=0, \quad s_{\Lambda_i}^1(\y)<s^*\}, 
  \label{gam1}\\
   \Gamma_{2,\Lambda_i}=\{(s_{\Lambda_i}^2,y) : \mathcal{G}_{\Lambda_i}(s_{\Lambda_i}^2(\y),y)=0, \quad s_{\Lambda_i}^2(\y)>s^*\}, 
  \label{gam1a}
   \end{align}
   where $s^*$ is given in Assumption \textbf{S-shaped} \eqref{assum2}.
   \end{lemma}

Proof:
Since $0<s<1$ and $\Lambda_i<0$ or $\Lambda_i>1$ the eigenvalues
$\lambda_{\Lambda_i}$ are well defined. Notice that $\frac{\partial \mathcal{G}_{\Lambda_i}(s,y)}{\partial s}=\frac{\partial ^2 f(s,y)}{\partial s ^2}+\frac{\lambda_{\Lambda_i}-\lambda_s}{s-\Lambda_i}$. Using the S-Shape assumption, we have $\frac{\partial \mathcal{G}_{\Lambda_i}(s,y)}{\partial s} \not= 0$, thus by the implicit function Theorem (IFT) there exist $s_{\Lambda_i}^1(\y)$ such that $\mathcal{G}_{\Lambda_i}(s_{\Lambda_i}^1(\y),y)=0$. Moreover, by the S-Shape assumption there exist $s_{\Lambda_i}^2(\y)$
satisfying Eq. $(\ref{gam1})$ of Lemma \ref{lema1}.$\Box$
 
 The coincidence between  two different $\lambda_{\Lambda_i}$
 occurs in some  special cases. 
 \begin{lemma}
		\label{comparacao}
 		At the states where
 		${\lambda}_{\Lambda_i}={\lambda}_{\Lambda_j}$, we have $f=s$ or $\Lambda_i=\Lambda_j$.
 	\end{lemma}

\subsection{Inflection loci}

Other structures appearing in this model are the inflection loci. 
The inflections are, generically, co-dimension 1 structures, in which the monotone increase of characteristic speed fails, i.e., $\nabla \lambda \cdot \vec{r}=0$.

In this model, we have $n-1$ fields. For the field $\lambda_s$, it is easy to see that $\nabla \lambda_s \cdot \vec{r}_s=\frac{1}{\varphi}\displaystyle{\frac{ \partial^2 \f}{\partial \s^2}}$, 
thus the inflection locus consists of the states  $s^*=s^*(\y)$, satisfying $\displaystyle{\frac{ \partial^2 \f}{\partial \s^2} (s^*,\cdot)=0}$. 
We denote the inflection for the field $(\lambda_s,\vec{r}_s)$ as $\mathcal{I}_s$:
\begin{equation}
   \mathcal{I}_s=\left\{(s^*,y)\text{ : } \frac{ \partial^2 \f}{\partial \s^2} (s^*,y)=0\right\}.
   \label{is}
\end{equation}
From IFT, if $\displaystyle{\frac{ \partial^3 \f}{\partial \s^3} (s^*,y)\neq 0}$ for $(s^*,y)\in \mathcal{I}_s$, then $\mathcal{I}_s$ is
a smooth hypersurface with co-dimension 1 in $\Omega$.
Moreover, if $\f$ does not depend on $y$, then $\mathcal{I}_s$ is a hyperplane in $\Omega$ with constant $s$. 

\begin{proposition}
 \label{inflectionproposition}
 The inflection locus $I_{\Lambda_i}$  of field $(\lambda_{\Lambda_i},\vec{r}_{\Lambda_i})$
 satisfies
 \begin{equation}
 \nabla \lambda_{\Lambda_i} \cdot \vec{r}_{\Lambda_i}=\frac{1}{\varphi}\frac{\s-\f}{(\s-\Lambda_i)^2}(\lambda_s-\gr{\lambda_{\Lambda_i}})\mathcal{H}_i(s,y),
\quad \text{ where } \quad \mathcal{H}_i =\vec{w_i}\cdot \vec{v_i},\label{fint}
 \end{equation}
with $\vec{w_i}=(\vec{w_i})_j$, for $j=1,\cdots,n-1$, and $(\vec{w_i})_j$  given by $(\ref{wee})$. Here $\vec{v_i}$ is a right eigenvector associated to generalized eigenvalue $\Lambda_i$ in \eqref{eigeqs2}. 
\end{proposition}
 	
 \textbf{Proof of Proposition $\ref{inflectionproposition}$}

Let $\lambda_{\Lambda_i}$ be given by $(\ref{eig2a})$. Calculating $\nabla \lambda_{\Lambda_i}$, we obtain
\begin{align}
 &\nabla \lambda_{\Lambda_i}=\frac{1}{\varphi}\left(\omega,
\varpi,\frac{\lambda_{\Lambda_i}}{u}\right), \text{ where } \notag \\
& \omega =\frac{1}{\s-\Lambda_i}\left(\left(\lambda_s-\lambda_{\Lambda_i}\right)+u\frac{\partial \Lambda_{i}}{\partial s}(\lambda_{\Lambda_i}-1)\right),\\
&\varpi_j=-\frac{u}{(s-\Lambda_i)^2}(f-s)\frac{\partial \Lambda_i}{\partial y_j}
+\frac{\partial f}{\partial y_j}\frac{1}{s-\Lambda_i}, \quad \text{ for } j=1,\cdots,n-1.
\end{align}
Using \eqref{e40} and $\vec{r}_{\Lambda_i}$ given in $(\ref{eiget1})$, we calculate $\nabla \lambda_{\Lambda_i} \cdot \vec{r}_{\Lambda_i}$. After  some calculations we obtain $(\ref{fint})$, 
for which $(\vec{w_i})_j$, for $j=1,\cdots,n-1$ are given by
\begin{align}
 (\vec{w_i})_j&=\omega (s-\Lambda_i)^2 a_j \delta_{1,j}+\frac{(s-\Lambda_i)^2\lambda_{\Lambda_i} b_j}{u}\delta_{n+1,j}+(s-\Lambda_i)^2\sum_{k=2}^{n}\delta_{k,j} \varpi_j.\label{wee}
\end{align} 
where $a_j$ and $b_j$ are given in \eqref{eA} and \eqref{eB}
 $\square$.  

Let us define the  hypersurfaces $\mathcal{C}_s$, $\Gamma_{\Lambda_i}$ and $\mathcal{J}_{\mathcal{H}_i}$ as the loci where $\s= \f$, $\lambda_s-\lambda_{\Lambda_i}=0$ and $\mathcal{H}_i(s,y)=0$, respectively. From Proposition \ref{inflectionproposition}
the inflection locus $\mathcal{I}_{\Lambda_i}$ is given by the union of above hypersurfaces. A similar decomposition of the inflection surface was obtained by \cite{helmut2011thermal}.

 \section{Discontinuous solutions}
\label{rsk}
Discontinuous solutions of systems of conservation laws
$(\ref{eq1})$ satisfy the Rankine-Hugoniot condition, i.e., for given left and right states
$(s^-,y^-,u^-)$ and $(s^+,y^+,u^+)$ we have
\begin{equation}
u^+F_i((s^+,y^+)-u^-F_i(s^-,y^-)=\sigma(G_i(s^+,y^+)-G_i(s^-,y^-)),
\label{rh1aa}
\end{equation}
with $i=1,\dots,n+1$, where  $\sigma$ is the propagation speed of the discontinuity.
Equation \eqref{rh1aa} can be rewritten as
\begin{align}
\Phi_{i}.[\sigma,u^+,u^-]=0, 
\label{rank1}
\end{align}
where $\Phi_{i}=(\Phi_{i1},\Phi_{i2},\Phi_{i3})$, with\\
$\Phi_{i1}=\s^+(\rho_{wi}^+-\rho_{oi}^+)+\rho_{oi}^+-(\s^-(\rho_{wi}^--\rho_{oi}^-)+\rho_{oi}^-)+(\rho_{ri}^+-\rho_{ri}^+)$,~\\
$\Phi_{i2}=-(\f^+(\rho_{wi}^+-\rho_{oi}^+)+\rho_{oi}^+)$,~\\ $\Phi_{i3}=(\f^-(\rho_{wi}^--\rho_{oi}^-)+\rho_{oi}^-)$ 
and  $\rho^+=\rho(\y^+)$, $\rho^-=\rho(\y^-)$,  $\f^+=\f(\s^+,y^+)$, $\f^-=\f(\s^-,y^-)$ with
$\rho=(\rho_{wi},\rho_{oi},\rho_{ri})$.

For each fixed state $(\s^-,y^-)$, the $\textit{Hugoniot locus}$ $\mathcal{HL}(\s^-,y^-)$ consists of all states $(\s^+,y^+)$
satisfying $(\ref{rank1})$. In \cite{lambert2019nonlinear1} it is proven that the Hugoniot locus satisfies
\begin{equation}
\mathcal{HL}(\s^-,y^-)=\left\{(\s^+,y^+): det(\Phi_{i}^T,\Phi_{k}^T,\Phi_{j}^T)=0\right\},
\label{hl1a}
\end{equation}
for all combinations of distinct indices $\{i,k,j\}\in \{1,\ldots,n+1\}$.
Also, we verify that instead of considering all index combination it is sufficient to reduce system
\eqref{hl1a} to $n-1$ linearly independent equations. This result can summarized as follows: let $i_1$ and $i_2$ $\in\{1,2,\cdots,n+1\}$ be two indices 
such that $\Phi_{i_1}$ and $\Phi_{i2}$ are linearly independent, then Equation \eqref{hl1a} reduces to
\begin{equation}
\mathcal{HL}(\s^-,y^-)=\left\{(\s^+,y^+): det(\Phi_{k}^T,\Phi_{i_1}^T,\Phi_{i_2}^T)=0\right\},
\label{hl2a1}
\end{equation}
for values of $k \in\{1,\ldots,n+1\}$ distinct from $i_1$ and $i_2$. The curve given by \eqref{hl2a1} is the intersection of $n-1$ independent surfaces in the
$n$ dimensional space $(s,y)$.

Assuming that there exist two indices $i,j$ with $i \neq j$ such that the denominator is different from zero for given states $(s^-,y^-)$ and $(s^+,y^+)$ then the values of $u^+$ can be obtained from \eqref{rh1aa} as
\begin{equation}
\frac{u^+}{u^-}=\frac{F_i(s^-,y^-)[G_j]-F_j(s^-,y^-)[G_i]}{F_i(s^+,y^+)[G_j]-F_j(s^+,y^+)[G_i]},
\label{key1}
\end{equation}
where $[G_i]=G_i(s^+,y^+)-G_i(s^-,y^-)$ and $[G_j]=G_j(s^+,y^+)-G_j(s^-,y^-)$. 

A particular branch of the Hugoniot locus can be characterized in the following

\begin{lemma}
For fixed $(s^-,y^-)$, a branch of $\mathcal{HL}(s^-,y^-)$ consists of the 
states of the form $(s,y)\in\Omega$, with $s$ variable and $y=y^-$. 
\end{lemma} 
 This branch is called \textit{Buckley-Leverett or saturation branch}.
\section{Elementary waves and interactions}

The  features of the coefficients $\rho_{rj}$ $(j=1,\cdots,n+1)$ in the accumulation function have implications on the quantitative behavior of the Riemann solutions. If such coefficients are constant in the phase space then the matrices $\mathcal{A}_1$ and $\mathcal{B}_1$ defined in \eqref{sd3} and \eqref{sd4} satisfy $\mathcal{A}_1\equiv0$ and $\mathcal{B}_1\equiv0$. This fact implies that
the matrix ${G}^{(2)}_{l,r}$ in \eqref{e40}
and the corresponding eigenvalues and eigenvectors of equation \eqref{eigeqs2} depend only on the variables $y$. In \cite{lambert2019nonlinear1} this case is studied and the existence of the splitting in the variables $u$ and $s$ from $y$ to
obtain the chemical rarefaction and the Hugoniot locus is shown. This splitting is extremely useful to obtain the wave curves.

In general to consider the coefficients $\rho_{rj}$ $(j=1,\cdots,n+1)$ as non constant functions in the accumulation terms in \eqref{acful} has significant implications in the structure of the solution. In such a case it is no longer possible to solve the Riemann problem in terms of the variables $y$ alone. Thus, considering surface complexes increases difficulties because more wave interactions appear in the construction of wave curves.

However, in general in the state variables $(s,y,u)$, the variable $u$ can be split from the variables $(s,y)$. Let us take $\Lambda$ as any of the eigenvalues $\Lambda_i$ ($i=1,\cdots,n-1$) given in \eqref{eigeqs2}. Using  the Proposition \ref{lemakt} and equation $(\ref{integral})$, we  obtain along the integral curves $\mathcal{W}^{\Lambda}$  the unknown variables $(s,y)$ independently of $u$ by solving
the system of differential equations:
\begin{equation}
\frac{ds}{d\xi}=r_{\Lambda}^1(s,y),\;
\frac{dy_i}{d\xi}=r_{\Lambda}^i(s,y),~ \text{for}~ i=2,\cdots,n,~~
\frac{du}{d\xi}=ug_{\Lambda}^{n+1}(s,y),\label{equvv}
\end{equation}
with $\xi=x/t$, $\xi^-=\lambda(s^-,y^-,u^-)$, $s(0)=s^-$, $y(0)=y^-$, $u(0)=u^-$  and $g_{\Lambda}^{n+1}(s,y)=r_{\Lambda}^{n+1}(s,y,u)/u$, where $r_{\Lambda}^{i}$ ($i=1,\cdots,n+1$) are given in \eqref{eiget1a}. Here the superscript index is used to indicate the coordinates of the vector $\vec{r}_\Lambda$ given in Eq. $(\ref{eig2a})$ of  Proposition \ref{lemakt}. After obtaining $(s(\xi),y(\xi))$, we use the expression from the  last
equation of $(\ref{equvv})$ yielding 
\begin{equation}
u=u^-exp({\gamma(\xi)}), ~~ \gamma(\xi)=\int_{\xi^-}^{\xi}{
g_{\Lambda}^{n+1}}(s(\eta),y(\eta))d\eta,\label{useses}
\end{equation}
where $\xi=x/t$, $\xi^-=\lambda(s^-,y^-)$ and $u^-$ is the initial
value of $u$ on the rarefaction wave, i.e. $u=u^-$ at $\xi
=\xi^-$. 
Furthermore, in \cite{lambert2009riemann}, the authors  proved that it is possible to project the Rankine-Hugoniot locus from a state $(s^-,y^-,u^-)$ 
in the $\Omega=\{(s,y)\}$ space. In this sense, we can solve the Riemann problem in the projection space $(s,y)$, which reduces the number of calculations.

From formula \eqref{useses}, we can deduce that if the coefficients of system \eqref{eq1} are such that the function $\gamma(\xi)$ assumes small values then the values of the variable $u$ is small. This case was observed in experiment studies, see  \cite{alvarez2018analytical1}.

We consider the interaction of the waves. There are two cases:
namely, waves in the same family, and waves in different families. The resonance condition leads to
the construction of curves formed by structures from distinct families. 

Here,
we have $n-1$ kinds of shocks, one associated to $\lambda_s$, which we denote by $\mathcal{H}_s$
and the other associated to $\lambda_{\Lambda_i}$ ($i=1,\cdots,n-1$), i.e., $\mathcal{H}_{\Lambda_i}$.
 Similarly to $\mathcal{R}_s$, the curve $\mathcal{H}_s$ is a straight
line parallel to the $s$ axis. On the other hand, for states in the neighborhood of $(s-, y-, u-)$, the chemical shocks $\mathcal{H}_{\Lambda_{i}}$ have
structures similar to the rarefaction $\mathcal{R}_{\Lambda_i}$. However, for states far from the left state, we do not identify the behavior of Rankine-Hugoniot locus. 

In summary, the elementary waves consist of
\begin{itemize}
    \item The B-L saturation wave  $\mathcal{R}_s$ where only the saturation varies. This wave curve stops
    at the inflection locus $I_s$ given in \eqref{is} and bifurcates at surface
    $\lambda_s=\lambda_{\Lambda_i}$ with $i=1,\ldots,n-1$,
    
    \item The shock wave curve $\mathcal{H}_s$ consists of the locus where only the saturation varies,
    
    \item The saturation wave curves $\mathcal{R}_{\Lambda_i}$ ($i=1,\ldots,n-1$) associated to the couple
    $(\lambda_{\Lambda_i},r_{\Lambda_i})$. These wave curves stop at the inflection loci $I_{\Lambda_i}$  given in \eqref{fint} and bifurcate at surfaces $\Gamma_{ij}$:
    $\lambda_{\Lambda_j}=\lambda_{\Lambda_i}$
    with $i,j=1,\ldots,n$ $i \neq j$ or  $\lambda_s=\lambda_{\Lambda_i}$ with $i=1,\ldots,n$,
    
    \item The composite wave curve $C_s=\mathcal{R}_s \cup \mathcal{H}_s$ consists of a characteristic B-L shock,
    
    \item The composite wave curves $C_{\Lambda_i}=\mathcal{R}_{\Lambda_i} \cup \mathcal{H}_{\Lambda_i}$, where $\mathcal{H}_{\Lambda_i}$ is a characteristic shock associated to the branch of the family $\lambda_{\Lambda_i}$ with $i=1,\ldots,n-1$,
    
    \item The constant state $\mathcal{C}$.
 
\end{itemize}

When the Riemann problem is solved the collision of waves of the same and different families may happen. Thus the possible collisions can be summarized as 

\begin{itemize}
\item The collision of two shocks: $\mathcal{H}_s \mathcal{H}_{\Lambda_i}$, $\mathcal{H}_{\Lambda_i}\mathcal{H}_{\Lambda_j}$ with $i \neq j$, $i,j=1,\cdots,n$ ,
\item The collision of two rarefaction waves: $\mathcal{R}_s \mathcal{R}_{\Lambda_i}$, $\mathcal{R}_{\Lambda_i}\mathcal{R}_{\Lambda_j}$ with $i \neq j$, $i,j=1,\cdots,n-1$,
\item The collision of rarefaction waves and a shock: $\mathcal{R}_s \mathcal{H}_{\Lambda_i}$, $\mathcal{R}_{\Lambda_i}\mathcal{H}_{\Lambda_j}$ with $i \neq j$, $i,j=1,\cdots,n-1$,
\item The collision of a shock and rarefaction waves: $\mathcal{H}_s \mathcal{R}_{\Lambda_i}$, $\mathcal{H}_{\Lambda_i}\mathcal{R}_{\Lambda_j}$ with $i \neq j$, $i,j=1,\cdots,n-1$.
\end{itemize}

Overtaking of elementary waves from the same family

\begin{itemize}
\item  i-shock overtakes another i-shock $\mathcal{H}_{\Lambda_i}\mathcal{H}_{\Lambda_i}$,
\item s-shock overtakes another s-shock $\mathcal{H}_s
\mathcal{H}_s$,
\item s-rarefaction wave overtakes i-shock $\mathcal{H}_{\Lambda_i} \mathcal{R}_s$,
\item s-shock overtakes i-rarefaction wave $\mathcal{R}_{\Lambda_i}\mathcal{H}_s$,
\item i-shock  wave overtakes i-rarefaction $\mathcal{R}_{\Lambda_i}\mathcal{H}_{\Lambda_i}$,
\item s-shock wave overtakes s-rarefaction $\mathcal{R}_{s}\mathcal{H}_{s}$,
\item i-shock overtakes i-rarefaction wave $\mathcal{R}_{\Lambda_i}\mathcal{H}_{\Lambda_i}$,
\item s-rarefaction wave overtakes another s-rarefaction wave $\mathcal{R}_{s}\mathcal{R}_{s}$.
\end{itemize}

\section{The Riemann Problem: Structure of waves}
\label{RP1}
The wave curves and bifurcation loci studied in the previous Sections are the main ingredients used to obtain the solution 
of $(\ref{eqt})$ for given initial data
\begin{equation}
\begin{cases}
 (s,y,u)_L & \text{ for } x<0,\\
 (s,y,\cdot)_R & \text{ for } x>0. 
\end{cases}
\end{equation}
For the system \ref{eqt} it is necessary to provide the initial 
data $u$ either on the right or left side. The missing value can be found from the other state and the wave sequences in the solution.

The Riemann solution consists of a sequence of elementary waves $w_k$ (shocks and  rarefactions) for $k=1,2,\cdots,m$
   and constant states $U_k$ for $k=1,2,\cdots,p$, in which $p$ is not known a priori.
   When the system is strictly hyperbolic and each
   field is genuinely nonlinear or linearly degenerate, for a system of $n$ equations, we can prove that we have a sequence of at most $n$ waves, see \cite{dafermos2005hyperbolic}.
   However, for the general case (with loss of strict hyperbolicity or loss of genuine nonlinearity) 
   it is not easy to determine a priori the number of waves in the solution. Here, we show that
   the system exhibits several nonlinear effects such as loss of hyperbolicity. 
   
 In any case, these waves are written as
\begin{equation}
\mathcal{U}_L\equiv \mathcal{U}_0 \overset{w_1}{\longrightarrow}
\mathcal{U}_1 \overset{w_2}{\longrightarrow}  \cdots
\overset{w_m}{\longrightarrow} \mathcal{U}_m\equiv \mathcal{U}_R
\label{sequence},
\end{equation}
where $\mathcal{U}=(s,y,u)$.   In
the Riemann solution it is necessary that the waves have increasing speed, the so called \textit{geometrical
compatibility}. Sometimes, this \textit{geometrical compatibility} 
is sufficient to furnish existence and uniqueness for the solution. Moreover, this condition is used to select the physical sequence of waves for the Riemann solution. 

With the previous considerations, we study the behavior of waves and
their interactions. In that direction, we study the Riemann solution for a particular problem, where we are able to obtain some interesting behavior of the wave curves, when these wave curves cross the inflection or coincidence loci. 

In the previous Sections we find bifurcations, eigenvalues and eigenvectors for system \eqref{eq1}. In summary, the eigenpairs for states $(s,\y)$, $\y=(\y_1,\cdots,\y_{n-1})$ are: a \textit{saturation} eigenpair of form
  $\displaystyle{\lambda_s=\frac{u}{\varphi}\frac{\partial f}{\partial s} \quad \vec{r}_s=(1,0,\cdots,0)}$
and $n-1$  compositional eigenvalues of form $
 \displaystyle{\lambda_{\Lambda}=\frac{u}{\varphi}\frac{f-\Lambda}{s-\Lambda}}$. We will also use the fact that each eigenvector generates its own rarefaction curves, which are denoted $\mathcal{R}_s$ for the field $\vec{r}_s$ 
 and $\mathcal{R}_{\Lambda_i}$ for the field $\vec{r}_{\Lambda_i}$.

For the rarefaction $\mathcal{R}_s$ only the saturation changes. The rarefactions associated to  $\mathcal{R}_{\Lambda_i}$ are in general transverse to $\mathcal{R}_s$. 

 In this model, we assume there is a coincidence surface between eigenvalues of different families, i.e., there is $(s,\y)$ for which
$\lambda_s=\lambda_{\Lambda_i}$.  This system loses   strict hyperbolicity on  this coincidence surface. 

\subsection{Application}
\label{hl2a}

In this section, we present numerical examples for a system of four conservation laws having the same structure as system \eqref{eqt}.
We give an overview of the main parameters of the model, a full description of which can be found in \cite{alvarez2018analytical1}. The physical model includes oil recovery improvement due to
the injection of carbonated water. The model also includes enhancement of the solubility of carbon dioxide in the water due to low salinity. We have four state variables. Two of these variables are the water saturation $s$ and the Darcy velocity $u$. 
The two other state variables are associated to the water phase, namely the hydrogen ion concentration $[H^{+}]$ (denoted here by $y_2$) and the chloride concentration $[Cl^-]$ (denoted here by $y_1$). 

The system has four state variables and four equations. These equations include the molar concentrations of the total hydrogen $(\rho_{w1})$, the total chloride $(\rho_{w2})$,  the total hydrogen minus twice the total oxygen $(\rho_{w3})$ with  carbon dioxide $(\rho_{o3})$ and decane $(\rho_{o4})$ (see \cite{alvarez2018analytical1}). They also include the two coefficients $\rho_{r3}$ and $\rho_{r4}$ in the accumulation term representing the formation of chemical complexes on the rock surface.

Disregarding the diffusive terms as in \cite{alvarez2018analytical1}, the accumulation  $G=(G_1,G_2,G_3,G_4)$ and flux functions $F=u\widehat{F}$  (with $\widehat{F}=(\widehat{F}_1,\widehat{F}_2,\widehat{F}_3,\widehat{F}_4)$)  are written as:
 \begin{align}
 \label{fsombre1}
 G&=( \varphi s \rho_{w1},\varphi s \rho_{w2},s \varphi \rho_{w3}  + 4\varphi(1-s)\rho_{o3}+(1-\varphi)\rho_{r3}, (1-s)\varphi\rho_{o4}+(1-\varphi)\rho_{r4}),
 \\
\widehat{F}
  &= (f\rho_{w1},f\rho_{w2},f\rho_{w3}+ 4(1-f)%
 \rho_{o3}, (1-f) \rho_{o4}).
 \label{fsombre}
 \end{align}
 Here, $\varphi$ is the porosity. The coefficient functions $\rho_{w1}$,
  $\rho_{w3}$,  $\rho_{o3}$, $\rho_{w2}$ and $\rho_{o4}$  are obtained through the PHREEQC program (see \cite{parkhurst2013user,parkhurst1999user}). We use the coefficients given in \cite{alvarez2018analytical1}. These coefficients depend on the  variables $y=(y_1,y_2)$ only. Here, we consider several possibilities for the coefficients $\rho_{r3}$ and $\rho_{r4}$. 
  
The fractional flow function $f$ has the $S$-shape with properties described in Assumption $\ref{assumption1}$. Here, we use the function presented in Section $\ref{fraction}$.

 We denote by $\lambda_k(A)$, with $k=s,\Lambda_1,\Lambda_2$  the characteristic speeds at state $A$. Also we denote by $\sigma(A,B)$ the shock speed from state $A$ to state $B$.
 
  In order to construct the Riemann solution in a region where $\lambda_{\Lambda_1} < \lambda_{s} < \lambda_{\Lambda_2}$, the geometrical compatibility condition requires that we start with a rarefaction wave $\mathcal{R}_{\Lambda_1}$, followed by either
 	a constant state or a shock wave (see below), which must
 	have a speed larger than the speed of the previous wave. On the other hand, if $\lambda_{\Lambda_2} > \lambda_{\Lambda_1} > \lambda_s$, we start with a saturation wave; either a rarefaction $\mathcal{R}_s$ or a shock $H_s$.  
 	
  	Finally, we can initiate building the solution
 	of the Riemann problem starting from the injection conditions of the reservoir (left state) until arriving at the initial conditions of the reservoir (right state). 
 	
 An in-house Riemann solver was developed to calculate and represent wave curves, which satisfy the compatibility and admissibility criteria (see \cite{alvarez2020resonance} for technical implementation details). In this solver, we automated the construction of slow and fast solution paths taking into account the bifurcation structures. Rarefaction curves are obtained by solving ODE systems, one for B-L rarefactions wave curves given in \eqref{ode1} and the other for chemical wave curves in \eqref{integral}. Moreover, shock waves are obtained by finding the intersection of the two surfaces constructed using \eqref{hl2a1} for $n=4$. Additional composition wave curves are necessary for crossing the inflection surfaces. When two eigenvalues coincide the rarefaction waves collide and numerical difficulties appear for the continuation of the wave curve. In this case the desingularizing algorithm developed in \cite{alvarez2020resonance} is used to allow wave curves to cross resonance surfaces.

   \subsubsection{Riemann solutions when the coefficients $\rho_{rj}$ are constant.}
   
 \textbf{\textit{Inflection and Coincidence Surfaces}}
 
\begin{figure}
      \begin{center}
      \includegraphics[scale=0.21]{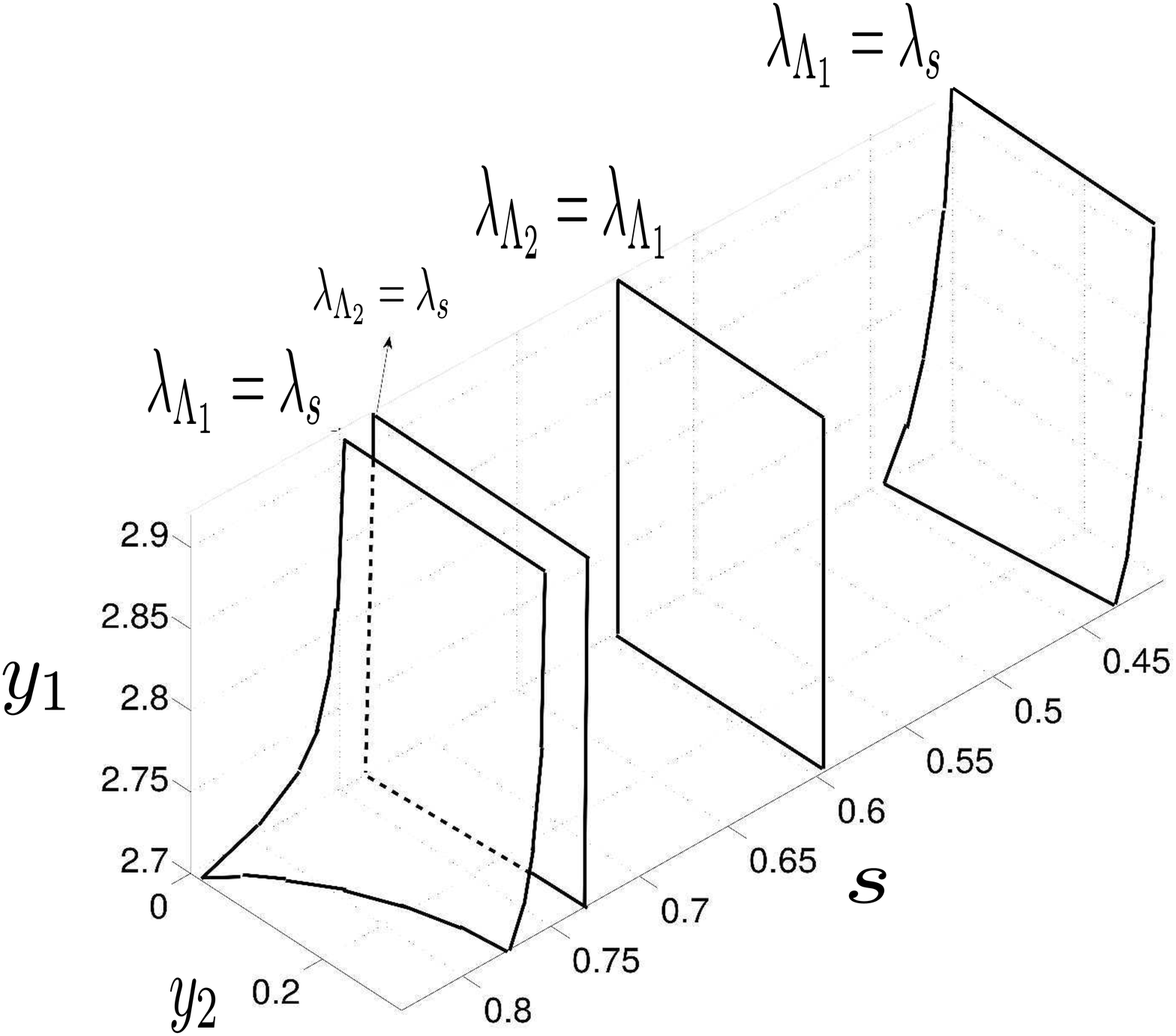}
       \caption{Coincidence Surfaces $\lambda_{\Lambda_i}=\lambda_s$, i=1,2 and $\lambda_{\Lambda_i}=\lambda_{\Lambda_j}$ for $i \neq j$. Each of these surfaces have almost constant saturation.}
       \label{surfaI1}
      \end{center}
      \end{figure}

 Rarefaction wave solutions suffer structural modifications along the loci where the eigenvalues coincide.  The coincidence surface has four connected parts given as the sets where
   speeds satisfy either $\lambda_{s}=\lambda_{\Lambda_1}$ (two parts) or  $\lambda_{\Lambda_2}=\lambda_{\Lambda_1}$ (one part) or $\lambda_{s}=\lambda_{\Lambda_2}$ (one part)  (see Figure \ref{surfaI1}). 

\begin{figure}
   \begin{center}
   \includegraphics[scale=0.21]{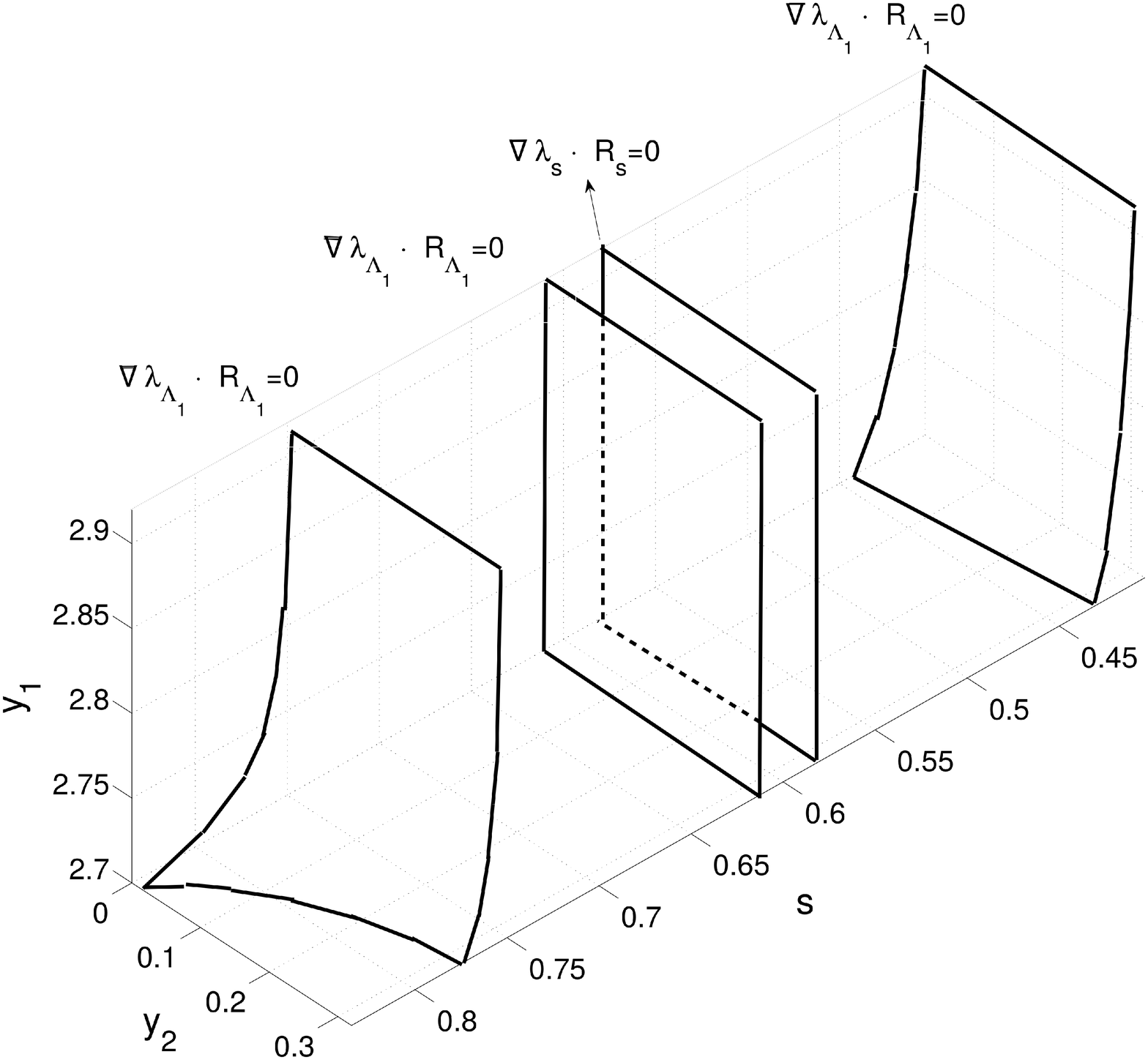}
   \caption{Inflection surfaces  for three eigenpairs $(\lambda_s,\vec{r}_s)$, $(\lambda_{\Lambda_1},\vec{r}_{\Lambda_1})$ and $(\lambda_{\Lambda_2},\vec{r}_{\Lambda_2})$.}
   \label{surfaI}
   \end{center}  
   \end{figure}

  From Proposition $\ref{inflectionproposition}$, we verify that the inflection surface consists of four connected parts: one part for  $\lambda_s$, two parts for   $\lambda_{\Lambda_1}$ and one part for $\lambda_{\Lambda_2}$ (see Figure \ref{surfaI}). 

  \textbf{\textit{Riemann Solutions}}
   
   In the general case, despite building elementary waves involves simple constructions, solving the Riemann problem is a difficult task due to the large number of combinations of elementary waves as well as to the presence of several bifurcation surfaces. However,
   when the coefficients $\rho_{jr}$, in the accumulation function in \eqref{fsombre1}, either are constant or have small derivatives relative to the other coefficients some useful simplifications occur.
  
   One of such simplifications is  that  the field associated to the eigenvalue $\lambda_{\Lambda_2}=({u}/{\varphi})({f}/{s})$ is linearly degenerate (see Corollary \ref{corolinfle}) and as a consequence the rarefaction and shock curves coincide. Moreover, curves associated to $\lambda_{\Lambda_1}$  and $\lambda_{s}$ 
  are invariant in the planes $y_1=(y_1)_C$ with constant $(y_1)_C$, while curves associated to $\lambda_{\Lambda_2}$ are transversal to these planes. This permits the displacement in the state space in the direction of the variable $y_1$.
 
     Let us to construct the solution for the Riemann-Goursat problem associated to the system of conservation laws $(\ref{eq1})$ with flux and accumulation functions given in \eqref{fsombre1}-\eqref{fsombre}, that is the solution of these equations with piecewise constant initial data
 \begin{equation}
   \left\{
   \begin{array}
   [c]{ll}%
   L=(s_{L},(y_1)_L,(y_2)_L,u_L) & \textbf{if}\hspace{0.2cm}x<0,\\
   R=(s_{R},(y_1)_R,(y_2)_R,\;\cdot\;) &
   \textbf{if}\hspace{0.2cm}x>0.
   \end{array}
   \right.\label{riemandata}
   \end{equation}
    To do that, we consider two examples taking into account the elementary waves and the bifurcation surfaces described above. In the first example the left state $L$ and the right state $R$ have different values for the variable $y_1$, i.e., $y_1=C_1$ and $y_1=C_2$, with constant $C_1$ and $C_2$. However, in the second example the left and right state have the same value for the variable  $y_1$. Because the variables $(s,y_1,y_2)$ can be split from the variable $u$, we are able to represent the Riemann solutions in the projected state space $(s,y_1,y_2)$ in $\{[0,1] \times [0,0.4] \times [2.7,4]\}$ (see \cite{lambert2019nonlinear1}). In the two examples we fixed $\rho_{3r}=0.01$ and $\rho_{4r}=0.02$.
  
  In both examples the left state $L$ belongs to the region where $\nabla \Lambda_s \cdot \vec{r}_s<0$ and $f>s$ and the right state $R$ belongs to the region where $\nabla \lambda_s \cdot \vec{r}_s>0$ and $f<s$. These restrictions imply that to cross the surface $\nabla \lambda_s \cdot \vec{r}_s=0$, we need only type B-L shocks. This is because $\nabla \lambda_s \cdot \vec{r}_s=0$ is a barrier for the B-L rarefactions.  Moreover, the surface $f=s$ constitutes a barrier for the remaining rarefaction curves (see \cite{lambert2019nonlinear1}). This fact can be verified by simple inspection of the formula for the infection surface given in \eqref{fint}.
  
 In the first example, we take the left state $L=(s_{L},(y_1)_L,(y_2)_L)=(0.99,0.03,2.74)$ belonging to a region where the condition $\lambda_{s} < \lambda_{\Lambda_1} < \lambda_{\Lambda_2}$ holds. We take the right state given by $R=(s_{R},(y_1)_R,(y_2)_R)=(0.27,0.3,2.74)$. In this case left and right states belong to  distinct planes i.e. $y_1=(y_1)_L=0.03$ and $y_1=(y_1)_R=0.3$.  
 
\begin{figure}
\begin{center}
\includegraphics[scale=0.22]{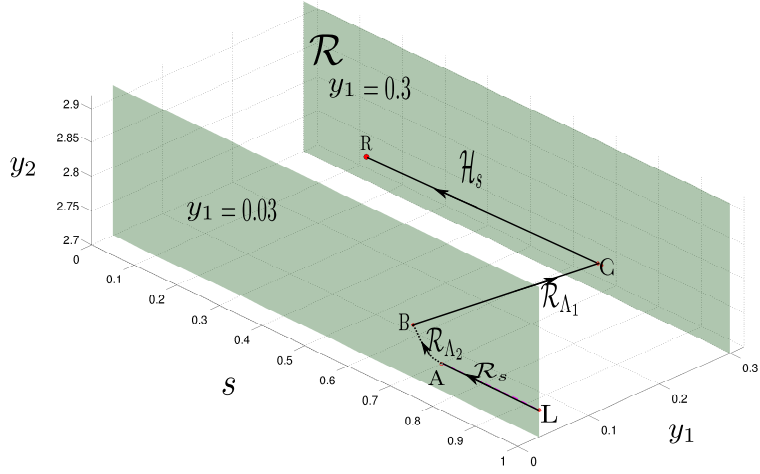}
\caption{Schematics of Riemann solution corresponding to left (L) and right states (R). The solution is obtained as a concatenation of four wave curves separated by two constant states: the first consists of a rarefaction saturation wave $\mathcal{R}_{s}(L)$ followed by a rarefaction curve $\mathcal{R}_{\Lambda_1}(A)$ (shortened $\mathcal{R}_{\Lambda_1}$ in the Figure), the second is a rarefaction curve $\mathcal{R}_{\Lambda_2}(B)$ (shortened $\mathcal{R}_{\Lambda_2}$ in the Figure) and the third is an admissible BL-shock $\mathcal{H}_s(C)$. The coordinates of each point are
$L=(0.99,0.03,2.74)$, $A=(0.77,0.03,2.74)$, $B=(0.70,0.03,2.77)$, $C=(0.70,0.3,2.74)$ and $R=(0.15,0.03,2.74)$. Denoting  $\lambda=(\Lambda_s,\lambda_{\Lambda_1},\lambda_{\Lambda_2})$, we have $\lambda(L)=(1.15 \times 10^{-5},2.66 \times 10^{-5},2.73 \times 10^{-5})$, $\lambda(A)=(2.14 \times 10^{-5},2.14 \times 10^{-5},3.36 \times 10^{-5})$, $\lambda(B)=(4.76 \times 10^{-5},2.47 \times 10^{-5},3.36 \times 10^{-5})$, $\lambda(C)=(5.26 \times 10^{-5},2.45 \times 10^{-5},3.42 \times 10^{-5})$, $\lambda(R)=(2.55 \times 10^{-6},2.58 \times 10^{-7},3.04 \times 10^{-5})$  and $\sigma(C,L)=5.26 \times 10^{-5})$.}
\label{satu2}
\end{center}
\end{figure}

 The solution of the Riemann problem can be obtained using the following wave curves. First, we take a saturation curve $\mathcal{R}_{s}(L)$ from $L$ to the state $A$ where $\lambda_{s}(A) =\lambda_{\Lambda_1}(A)$ (See Figure \ref{satu2}). From the state $A$, we take a rarefaction curve $\mathcal{R}_{\Lambda_1}(A)$ connected to the physical boundary. 
 
\begin{figure}
\begin{center}
\includegraphics[scale=0.39]{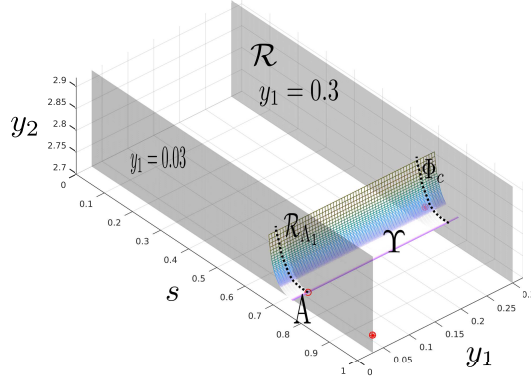}
\caption{Here we show the auxiliary surface $\Upsilon$ made up of rarefaction curves of type $\mathcal{R}_{\Lambda_2}(B^*)$ starting at all points $B^*$ on the rarefaction $\mathcal{R}_{\Lambda_1}(A)$ (shortened $\mathcal{R}_{\Lambda_1}$ in the Figure). The plane $\mathcal{R}$ corresponds to $y_2=0.3$. Here $\Phi_c$ represents the intersection curve between the surface $\Upsilon$ and the plane $\mathcal{R}$.}
\label{satu2aa}
\end{center}
\end{figure}

 From here the continuation of curve $\mathcal{R}_{\Lambda_1}(A)$ requires the determination of a point belonging to both curves $\mathcal{R}_{\Lambda_2}(B^*)$ and $\mathcal{H}_s(R)$ defined below. Since seeking such a point in the three-dimensional space is a numerically  unstable procedure, the construction of an auxiliary surface $\varUpsilon$ is necessary (see Figure \ref{satu2aa}). Thus in this example, we construct a surface $\Upsilon$ made up of rarefaction curves  $\mathcal{R}_{\Lambda_2}(B^*)$ starting at all points $B^*$
 on the rarefaction $\mathcal{R}_{\Lambda_1}(A)$.
 In this way, such curves form a surface which connects the rarefaction curve $\mathcal{R}_{\Lambda_1}(A)$ with the plane $\mathcal{R}$ ($y_1=0.3$), determining the intersection curve $\Phi_c$.
Now, we take a backward B-L shock curve out of the right
state $R=(0.27,0.3,2.74)$ on which the characteristic speed decreases. Afterwards, we look for the intersection state $C$ of the auxiliary surface $\Upsilon$ with the
backward shock curve $\mathcal{H}_s(R)$ from $R$ to $C$ (see Figure \ref{satu2}).
The existence of the state $C$ is guaranteed by the facts that 
in the curve $\mathcal{H}_s(R)$ only the saturation variable changes and
the auxiliary surface $\Upsilon$ has an intersection curve $\Phi_c$ with the plane $\mathcal{R}$ ($y_1=0.3$), which passes through the point with the same coordinate $y_2=2.74$ of the right state $R$.

Now, we look for the rarefaction $\mathcal{R}_{\Lambda_2}(B)$ on the surface $\Upsilon$ from the point $B$ to $C$ (See Figure \ref{satu2}), where the state $B$ is chosen such that the curve $\mathcal{R}_{\Lambda_2}(B)$ crosses the plane $\mathcal{R}$, $y_1=0.3$ at the point $C$ with coordinate $y_2=2.74$. Here, we use the fact that the curve $\mathcal{R}_{\Lambda_2}(B)$ is a contact wave, therefore the characteristic speed along this curve remains constant equal to $\lambda_{\Lambda_1}(B)$. This fact guarantees compliance to the monotone increase of the wave speed.

Since at the state $C$ the wave speed satisfies $\lambda_{\Lambda_2}(C) = \sigma(C,R)$, we have at this point a characteristic BL-shock $\mathcal{H}_s(C)$ from state $C$ to $R$. Finally, the solution is given by the wave sequence
$L\xrightarrow{\mathcal{R}_{s}}(A)\xrightarrow{\mathcal{R}_{\Lambda_1}}(B)
 \xrightarrow{\mathcal{R}_{\Lambda_2}} (C) \xrightarrow{\mathcal{H}_s}(R)$.
 
Similarly, it can be verified that the wave sequences $L\xrightarrow{\mathcal{R}_{s}}(A)\xrightarrow{\mathcal{R}_{\Lambda_1}}(B)
 \xrightarrow{\mathcal{H}_s}(R^*)$ solve the Riemann problem
 from the left state $L$ to all right states $R^*$ given by $R^{*}=(s,0.3,2.74)$ where $s \in (0.15,0.26)$ and the shock speeds satisfy $\sigma(C,R^*) \in (4.13 \times 10^{-5}, 5.11 \times 10^{-5})$. We verify that $\sigma(C,R^{*}) > \lambda_{\Lambda_2}(C)$ holds for all states $R^*$. Moreover, the shocks from state $C$ to $R^*$ are 3-Lax shocks (see \cite{lax1957hyperbolic,lax1973hyperbolic}). We verify that at each intermediate state of the wave curves the geometrical compatibility condition holds. We also check the possibility of the existence of other combinations of curves from the left to right states, but these cases do not satisfy such condition. However, we do not have a theoretical result which guarantees that the solution found constitutes the unique solution for this problem for chosen $L$ and $R$. 

\begin{figure}
\begin{center}
\includegraphics[scale=0.22]{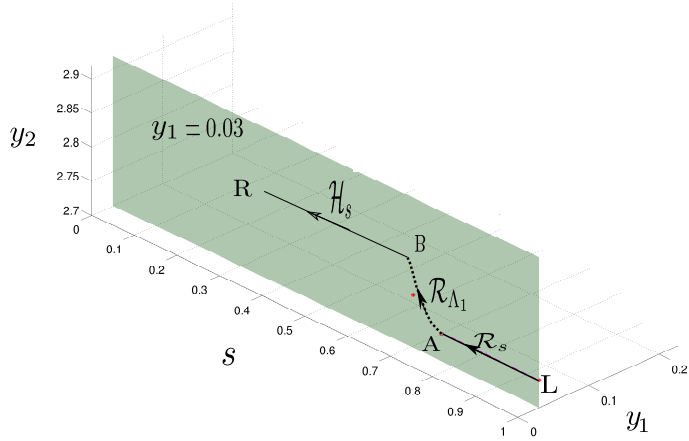}
\caption{Schematics of Riemann solution corresponding to left (L) and right states (R). The solution is obtained from a concatenation of three wave curves separated by two constant states: the first consists of a rarefaction saturation curve $R_{s}(L)$ (shortened $\mathcal{R}_{s}$ in the Figure) followed by a compatible rarefaction curve $\mathcal{R}_{\Lambda_1}(A)$ (shortened $\mathcal{R}_{\Lambda_1}$ in the Figure), the second one is the rarefaction curve $\mathcal{R}_{\Lambda_1}(B)$ and the third one is an admissible BL-shock $\mathcal{H}_s(C)$. The coordinates of each point are
$L=(0.99,0.03,2.74)$, $A=(0.77,0.03,2.74)$, $B=(0.68,0.03,2.85)$ and $R=(0.33,0.03,2.85)$. Let us denote  $\lambda=(\lambda_s,\lambda_{\Lambda_1},\lambda_{\Lambda_2})$, we have $\lambda(L)=(1.15 \times 10^{-5},2.66 \times 10^{-5},2.73.15 \times 10^{-5})$, $\lambda(A)=(2.14 \times 10^{-5},2.14 \times 10^{-5},3.36 \times 10^{-5})$, $\lambda(B)=(6.01 \times 10^{-5},2.70 \times 10^{-5},3.28 \times 10^{-5})$, $\lambda(R)=(7.96 \times 10^{-6},1.02 \times 10^{-6},2.7 \times 10^{-5})$  and $\sigma(C,L)=6.01 \times 10^{-5})$.}
\label{satu3}
\end{center}
\end{figure}

In the second example, we take the left state $L=(s_{L},(y_1)_L,(y_2)_L)=(0.99,0.03,2.74)$ satisfying  $\lambda_{s}(L) < \lambda_{\Lambda_1}(L) < \lambda_{\Lambda_2}(L)$. Assuming that the right state is given by $R=(s_{R},(y_1)_R,(y_2)_R)=(0.33,0.03,2.85)$, with $L$ and $R$ belonging to the same plane $\{y_1=0.03\}$,  the solution of the Riemann problem can be obtained with the following sequence of curves.
First, from $L$ to the state $A$ where $\lambda_{s}(A) = \lambda_{\Lambda_1}(A)$, we take the saturation wave $R_{s}(L)$  (see Figure \ref{satu3}). From the state $A$, we take the rarefaction $\mathcal{R}_{\Lambda_1}(A)$
connecting to the physical boundary.

To  reach the right state $R$ from here, it is useful to construct an extension of the rarefaction $R_{\Lambda_1}(A)$. Such extension curve is defined as those points $X$ belonging to a shock between states $X$ and $B$ in $R_{\Lambda_1}(A)$ such that $\sigma(X,B)=\lambda_s(B)$.

In this way, we find an intermediate state $B$ such that $\sigma(B,R)=\lambda_s(B)$
which consists of a characteristic B-L shock to the point $R$.
Finally, the solution is given by the wave sequence 
$L\xrightarrow{\mathcal{R}_s}(A)\xrightarrow{\mathcal{R}_{\Lambda_1}}(B)
 \xrightarrow{\mathcal{H}_s}(R)$. We verify that at the intermediate states $A$ and
 $B$ the compatibility condition of increasing speed is satisfied.

Other calculations were made using nonconstant coefficients $\rho_{jr}$ depending on the variables $y_1$ and $y_2$. In such examples changes of the intermediate states are observed but there are no changes with the respect to the constant coefficients case in the structure of the wave curves. However, these calculations are not sufficient to elucidate the issue of structural stability of the solutions.
This issue will be studied in future work.

\section{Conclusion}
\label{con}

In this paper the eigenvalue analysis of a system of $n+1$ conservation laws is performed. A method to build the discontinuous solutions and the bifurcation surfaces is presented as well. These elementary waves and surfaces are the basis for the subsequent analysis of the solution paths. It is noteworthy that these solution paths are useful for verification of numerical methods.

An overview of possible wave collisions is presented.
Using the elementary waves and the bifurcation surfaces some Riemann solutions are found. Due to the complexity of the problem in general, some hypotheses on the coefficients are made to obtain algorithms that allow easy application  of the wave curve method.  Structural stability of the solution of this kind of system is an open question.

\section*{Acknowledgement}
The authors thank Sergio Pilotto as well as support from CAPES  under grant 88881.156518/2017-01, CAPES/NUFFIC grant 88887.156517/2017-00, CNPq under grants
405366/2021-3, 306566/2019-2, FAPERJ under grants E-
26/210.738/2014, E-26/202.764/2017, E-26/201.159/2021.

\bibliographystyle{plain}
\bibliography{Surfacecomplexx1}

\end{document}